 \documentclass[12pt]{article}

 \usepackage{mathrsfs,amsfonts}
 \usepackage[dvips]{color}

 \newtheorem{theorem}{Theorem}[section]
 \newtheorem{lemma}[theorem]{Lemma}
 \newtheorem{corollary}[theorem]{Corollary}
 \newtheorem{proposition}[theorem]{Proposition}
 \newtheorem{example}[theorem]{Example}

 \def\blemma{\begin{lemma}\sl{}\def\elemma{\end{lemma}}}
 \def\bproposition{\begin{proposition}\sl{}\def\eproposition{\end{proposition}}}
 \def\btheorem{\begin{theorem}\sl{}\def\etheorem{\end{theorem}}}
 \def\bcorollary{\begin{corollary}\sl{}\def\ecorollary{\end{corollary}}}

 \def\beqlb{\begin{eqnarray}}\def\eeqlb{\end{eqnarray}}
 \def\beqnn{\begin{eqnarray*}}\def\eeqnn{\end{eqnarray*}}

 \def\proof{\noindent{\it Proof.~~}}\def\qed{\hfill$\Box$\medskip}

 \def\<{\langle}\def\>{\rangle}

 \def\mcr{\mathscr}\def\mbb{\mathbb}\def\mbf{\mathbf}

 \def\ar{\!\!&}\def\nnm{\nonumber}

 \def\lline{---------------------------------}

\begin{document}

\ %\noindent{08nonneg (draft form: 2008/01/26)}

%\tableofcontents

\bigskip\bigskip

\centerline{\Large\bf Stochastic equations of non-negative}

\smallskip

\centerline{\Large\bf processes with jumps\footnote{Supported by
NSFC (10525103 and 10721091).}}

\bigskip\bigskip

\centerline{Zongfei Fu}

\centerline{School of Information Sciences, Renmin University of
China,}

\centerline{Beijing 100872, People's Republic of China}

\centerline{E-mail: \tt fuzf@ruc.edu.cn}

\medskip

\centerline{Zenghu Li}

\centerline{School of Mathematical Sciences, Beijing Normal
University,}

\centerline{Beijing 100875, People's Republic of China}

\centerline{E-mail: \tt lizh@bnu.edu.cn}

\bigskip

\centerline{\lline\lline\lline}

\medskip

{\narrower{\narrower

\noindent\textit{Abstract.} We study stochastic equations of
non-negative processes with jumps. The existence and uniqueness of
strong solutions are established under Lipschitz and non-Lipschitz
conditions. The comparison property of two solutions are proved
under suitable conditions. The results are applied to stochastic
equations driven by one-sided L\'evy processes and those of
continuous state branching processes with immigration.

\smallskip


\smallskip


\smallskip

\noindent\textbf{Abbreviated Title:} Stochastic equations of
non-negative processes

\par}\par}

\bigskip

\centerline{\lline\lline\lline}

\bigskip

%%%%%%%%%%%%%%%%%%%%%%%%%%%%%%%%%%%%%%%

\section{Introduction}

\setcounter{equation}{0}

Stochastic differential equations with jumps have been playing
increasingly important roles in various applications. Under
Lipschitz conditions, the existence and uniqueness of strong
solutions of jump-type stochastic equations can be established by
arguments based on Gronwall's inequality and the results on
continuous type equations; see e.g.\ Ikeda and Watanabe (1989). In
view of the result of Yamada and Watanabe (1971), weaker conditions
would be sufficient for the existence and uniqueness of strong
solutions for one-dimensional equations. As an example of jump-type
equations, let us consider the simple equation
 \beqlb\label{1.1}
dx(t) = \phi(x(t-))dz(t), \qquad t\ge0.
 \eeqlb
By a result of Bass (2003), if $\{z(t)\}$ is a symmetric stable
process of order $1<\alpha<2$ and if $\phi(\cdot)$ is a bounded
continuous function with modulus of continuity $\rho(\cdot)$
satisfying
 \beqlb\label{1.2}
\int_{0+} \frac{1}{\rho(x)^\alpha}dx = \infty,
 \eeqlb
then (\ref{1.1}) admits a strong solution and the solution is
pathwise unique. This condition is exactly the analogue of the
Yamada-Watanabe criterion for the diffusion coefficient. When the
integral in (\ref{1.2}) is finite, Bass (2003) constructed a
continuous function $\phi(\cdot)$ having continuity modulus
$\rho(\cdot)$ for which the pathwise uniqueness for (\ref{1.1})
fails; see also Bass \textit{et al.}\ (2004). In view of
(\ref{1.2}), if a power function $\rho(x) = \mbox{const}\cdot
x^\beta$ applies for all symmetric stable processes with parameters
$1<\alpha<2$, we must have $\beta=1$. In other words, a universal
continuity modulus condition for jump-type stochastic differential
equations would not be a great improvement of the Lipschitz
condition.

For equations driven by non-symmetric noises, there is a new
difficulty brought about by the compensators of the noises. For
example, let us consider the equation (\ref{1.1}) again with
$\{z(t)\}$ being a one sided stable process of order $1<\alpha<2$.
For any $\varepsilon>0$, let
 \beqnn
z_\varepsilon(t) = \sum_{0< s\le t}[z(s) - z(s-)] 1_{\{z(s) - z(s-)>
\varepsilon\}}
 \eeqnn
and let $c_\varepsilon = \mbf{E}[z_\varepsilon(1)]$. We can define
another centered L\'evy process $\{w_\varepsilon(t)\}$ by
 \beqnn
w_\varepsilon(t) = z(t) - z_\varepsilon(t) + c_\varepsilon t.
 \eeqnn
Between any two neighboring jumps of $\{z_\varepsilon(t)\}$,
equation (\ref{1.1}) reduces to
 \beqlb\label{1.3}
dx(t) = \phi(x(t-))dw_\varepsilon(t) - c_\varepsilon\phi(x(t-))dt.
 \eeqlb
Then one would expect that, in order that the pathwise uniqueness
holds for (\ref{1.1}) or (\ref{1.3}), the function $\phi(\cdot)$
should be at least as regular as the drift coefficient in the
Yamada-Watanabe criterion. In other words, it should possess a
continuity modulus $r(\cdot)$ satisfying
 \beqlb\label{1.4}
\int_{0+} \frac{1}{r(x)}dx
 =
\infty,
 \eeqlb
which is much stronger than (\ref{1.2}).

Continuous state branching processes with immigration
(CBI-processes) constitute an important class of non-negative Markov
process with non-negative jumps. They were introduced in Kawazu and
Watanabe (1971) as approximations of classical Galton-Watson
branching processes with immigration. Many interesting applications
of them have been found since then. In particular, CBI-processes are
known as Cox-Ingersoll-Ross models (CIR-models) and have been used
widely in the study of mathematical finance; see, e.g., Duffie
\textit{et al.}\ (2003) and Lamberton and Lapeyre (1996). Up to a
minor moment assumption, a conservative CBI-process has generator
$A$ defined by
 \beqlb\label{1.5}
Af(x)
 \ar=\ar
axf^{\prime\prime}(x) + \int_0^\infty \big[f(x+z) - f(x)
- zf^\prime(x)\big]x \nu_0(dz) \nnm \\
 \ar \ar
+\, (b+\beta x)f^\prime(x) + \int_0^\infty \big[f(x+z)-f(x)\big]
\nu_1(dz),
 \eeqlb
where $a\ge 0$, $b\ge 0$ and $\beta$ are constants, and $\nu_0(dz)$
and $\nu_1(dz)$ are $\sigma$-finite measures on $(0,\infty)$
satisfying
 \beqlb\label{1.6}
\int_0^\infty (z\land z^2) \nu_0(dz) + \int_0^\infty (1\land z)
\nu_1(dz)< \infty.
 \eeqlb
Let $\{B(t)\}$ be a standard Brownian motion and let
$\{N_0(ds,dz,du)\}$ and $\{N_1(ds,dz)\}$ be Poisson random measures
with intensities $ds\nu_0(dz)du$ and $ds\nu_1(dz)$, respectively.
Suppose that $\{B(t)\}$, $\{N_0(ds,dz,du)\}$ and $\{N_1(ds,dz)\}$
are independent of each other. Let $\tilde{N}_0(ds,dz,du)$ be the
compensated measure of $N_0(ds,dz,du)$. From Theorems~5.1 and~5.2 of
Dawson and Li (2006) it follows that the stochastic equation
 \beqlb\label{1.7}
x(t)
 \ar=\ar
x(0) + \int_0^t \sqrt{2ax(s)}dB(s) + \int_0^t\int_0^\infty
\int_0^{x(s-)} z \tilde{N}_0(ds,dz,du) \nnm \\
 \ar \ar
+ \int_0^t (b+\beta x(s))ds + \int_0^t\int_0^\infty z N_1(ds,dz);
 \eeqlb
a unique non-negative strong solution; see also Corollary~\ref{c6.2}
of this paper. By It\^o's formula it is easy to see the solution
$\{x(t)\}$ of (\ref{1.7}) is a CBI-process with generator given by
(\ref{1.5}).

The purpose of the present paper is to study stochastic equations of
non-negative processes with jumps that generalizes the CBI-processes
described above. This kind of processes arise naturally in
applications, but they are excluded by most of the existing results
in the literature because of the degeneracy of their generators. By
specifying to non-negative processes we can make the best use of the
first moment analysis, which turns out to be essential for a number
of applications. We provide here some criteria for the existence and
uniqueness of strong solutions of those equations. The main idea of
those criteria is to assume some monotonicity condition on the
kernel associated with the compensated noise so that the continuity
conditions can be released. It follows from our criterion that
(\ref{1.1}) has a unique non-negative strong solution if $\{z(t)\}$
is a one-sided $\alpha$-stable process and if $\phi(\cdot)$ is a
$(1/2)$-H\"older continuous non-decreasing function satisfying
$\phi(0)=0$.

To describe another consequence of our criteria, suppose that
$1<\alpha<2$ and $(a,b,\beta,\nu_1)$ are given as above. Let
$\{B(t)\}$ be a standard Brownian motion, $\{z_0(t)\}$ be a
one-sided $\alpha$-stable process with characteristic measure
$z^{-1-\alpha}dz$ and $\{z_1(t)\}$ be a non-decreasing L\'evy
process with characteristic measure $\nu_1(dz)$. Suppose that
$\{B(t)\}$, $\{z_0(t)\}$ and $\{z_1(t)\}$ are independent of each
other. We shall see that for any $c\ge 0$ there is a unique
non-negative strong solution to
 \beqlb\label{1.8}
dx(t) = \sqrt{2ax(t)} dB(t) + \sqrt[\alpha]{cx(t-)} dz_0(t) + (\beta
x(t)+b)dt + dz_1(t).
 \eeqlb
A particular case of this equation has been considered in Lambert
(2007), where the uniqueness of the solution was left open. The
solution $\{x(t)\}$ of (\ref{1.8}) is a CBI-process with generator
given by (\ref{1.5}) with $\nu_0(dz)$ replaced by $c z^{-1-\alpha}
dz$. Of course, for any $c>0$ and $1<\alpha<2$ the continuity
modulus of the coefficient $x \mapsto \sqrt[\alpha]{cx}$ in
(\ref{1.8}) does not satisfy condition (\ref{1.4}).

Some basic results for stochastic equations of non-negative
processes are provided in Section~2. In particular, we give a
Lipschitz condition for the existence and uniqueness of the strong
solution. In Section~3 the pathwise uniqueness is studied. In
Section~4 we prove a weak existence result by second moment
arguments. The existence and uniqueness of strong solutions under
non-Lipschitz conditions are established in Section~5, where only
first moment conditions are required. We also prove two simple
properties of the solutions, continuous dependence on the initial
value and comparison property. In Section~6, we illustrate some
applications of the main results (Theorems~\ref{t2.5}, \ref{t5.1}
and~\ref{t5.2}).

Throughout this paper, we assume $(\Omega, \mcr{F}, \mcr{F}_t,
\mbf{P})$ is a filtered probability space satisfying the usual
hypotheses. Moreover, we make the conventions
 \beqnn
\int_a^b = \int_{(a,b]}
 \quad\mbox{and}\quad
\int_a^\infty = \int_{(a,\infty)}
 \eeqnn
for any real numbers $a\le b$.

The theory of jump-type stochastic equations is not as well
developed as that of continuous ones; see Bass (2004), Jacob and
Schilling (2001) and the references therein. We refer to Ikeda and
Watanabe (1989) and Protter (2003) for the theory of stochastic
analysis and to Bertoin (1996) and Sato (1999) for the theory of
L\'evy processes.

%%%%%%%%%%%%%%%%%%%%%%%%%%%%%%%%%%%%%%%%%%%%%%

\section{Preliminaries and Lipschitz conditions}

\setcounter{equation}{0}

In this section, we prove some basic results on stochastic equations
of non-negative processes with jumps. Let $U_0$ and $U_1$ be
complete separable metric spaces, and let $\mu_0(du)$ and
$\mu_1(du)$ be $\sigma$-finite Borel measures on $U_0$ and $U_1$,
respectively. Suppose that
 \begin{itemize}

\item
$x\mapsto \sigma(x)$ and $x\mapsto b(x)$ are continuous functions on
$\mbb{R}$ satisfying $\sigma(x) =0$ and $b(x) \ge 0$ for $x\le 0$;

\item
$(x,u) \mapsto g_0(x,u)$ is a Borel function on $\mbb{R} \times U_0$
such that $g_0(x,u)+x\ge 0$ for $x>0$, and $g_0(x,u)=0$ for $x\le
0$;

\item
$(x,u) \mapsto g_1(x,u)$ is a Borel function on $\mbb{R} \times U_1$
such that $g_1(x,u)+x\ge 0$.

 \end{itemize}
Let $\{B(t)\}$ be a standard $(\mcr{F}_t)$-Brownian motion and let
$\{p_0(t)\}$ and $\{p_1(t)\}$ be $(\mcr{F}_t)$-Poisson point
processes on $U_0$ and $U_1$ with characteristic measures
$\mu_0(du)$ and $\mu_1(du)$, respectively. Suppose that $\{B(t)\}$,
$\{p_0(t)\}$ and $\{p_1(t)\}$ are independent of each other. Let
$N_0(ds,du)$ and $N_1(ds,du)$ be the Poisson random measures
associated with $\{p_0(t)\}$ and $\{p_1(t)\}$, respectively. Let
$\tilde{N}_0(ds,du)$ be the compensated measure of $N_0(ds,du)$. By
a \textit{solution of} the stochastic equation
 \beqlb\label{2.1}
x(t)
 \ar=\ar
x(0) + \int_0^t \sigma(x(s-))dB(s) + \int_0^t\int_{U_0} g_0(x(s-),u)
\tilde{N}_0(ds,du) \nnm \\
 \ar \ar
+ \int_0^t b(x(s-))ds + \int_0^t\int_{U_1} g_1(x(s-),u) N_1(ds,du)
 \eeqlb
we mean a c\`adl\`ag and $({\mcr{F}}_t)$-adapted process $\{x(t)\}$
that satisfies the equation almost surely for every $t\ge0$. Since
$x(s-) \neq x(s)$ for at most countably may $s\ge 0$, we can also
use $x(s)$ instead of $x(s-)$ for the integrals with respect to
$dB(s)$ and $ds$ on the right hand side of (\ref{2.1}).

\bproposition\label{p2.1} If $\{x(t)\}$ satisfies (\ref{2.1}) and
$\mbf{P}\{x(0)\ge 0\} =1$, then $\mbf{P}\{x(t)\ge 0$ for all $t\ge
0\} =1$. \eproposition

\proof Suppose there exists a constant $\varepsilon>0$ so that $\tau
:= \inf\{t\ge 0: x(t)\le -\varepsilon\}< \infty$ with non-zero
probability. It is easy to see that $x(\tau-) = x(\tau) = -
\varepsilon$ on the event $\{\tau< \infty\}$. Let $\sigma = \inf\{s<
\tau: x(t)\le 0$ for all $s\le t\le \tau\}$. Then $\sigma< \tau$, so
we can choose a deterministic time $r\ge 0$ so that $\{\sigma\le r<
\tau\}$ has non-zero probability. On the event $\{\sigma\le r<
\tau\}$ we have
 \beqnn
x(t\land \tau) = x(r\land \tau) + \int_{r\land \tau}^{t\land \tau}
b(0)ds + \int_{r\land \tau}^{t\land \tau} \int_{U_1} g_1(0,u)
N_1(ds,du), \qquad t\ge r,
 \eeqnn
so $t\mapsto x(t\land \tau)$ is non-decreasing. Since $x(r)>
-\varepsilon$ on $\{r< \tau\}$, we get a contradiction. \qed

In the sequel, we shall always assume the initial variable $x(0)$ is
non-negative, so Proposition~\ref{p2.1} implies that any solution of
(\ref{2.1}) is non-negative. Then we can assume the ingredients are
defined only for $x\ge 0$. For non-negative processes we can use the
first moment estimates, which is essential for the CBI-process and
solutions of stochastic equations of driven by one-sided stable
noises. For the convenience of the statements of the results, let us
introduce the following conditions:
 \begin{itemize}

\item[{\rm(2.a)}]
There is a constant $K\ge 0$ and a subset $U_2\subset U_1$ such that
$\mu_1(U_1\setminus U_2)< \infty$ and
 \beqnn
|b(x)| + \int_{U_2} \sup_{0\le y\le x} |g_1(y,u)| \mu_1(du) \le
K(1+x), \qquad x\ge 0;
 \eeqnn

\item[{\rm(2.b)}]
There is a non-negative and non-decreasing function $x\mapsto L(x)$
on $\mbb{R}_+$ so that
 \beqnn
\sigma(x)^2 + \int_{U_0} \sup_{0\le y\le x}[|g_0(y,u)|\land
g_0(y,u)^2] \mu_0(du) \le L(x), \qquad x\ge 0.
 \eeqnn

 \end{itemize}

\bproposition\label{p2.2} Suppose that conditions (2.a,b) are
satisfied. If there is a strong solution to
 \beqlb\label{2.2}
x(t)
 \ar=\ar
x(0) + \int_0^t \sigma(x(s))dB(s) + \int_0^t\int_{U_0}
g_0(x(s-),u) \tilde{N}_0(ds,du) \nnm \\
 \ar \ar
+ \int_0^t b(x(s))ds + \int_0^t\int_{U_2} g_1(x(s-),u) N_1(ds,du),
 \eeqlb
there is also a strong solution to (\ref{2.1}). If the pathwise
uniqueness of solution holds for (\ref{2.2}), it also holds for
(\ref{2.1}). \eproposition

\proof The results hold trivially if $\mu_1(U_1\setminus U_2)= 0$,
so we assume $0< \mu_1(U_1\setminus U_2)< \infty$ in this proof.
Suppose that (\ref{2.2}) has a strong solution $\{x_0(t)\}$. Let
$\{S_k: k=1,2,\cdots\}$ be the set of jump times of the Poisson
process
 \beqnn
t \mapsto \int_0^t\int_{U_1\setminus U_2} N_1(ds,du).
 \eeqnn
We have clearly $S_k\to \infty$ as $k\to \infty$. For $0\le t<S_1$
set $y(t) = x_0(t)$. Suppose that $y(t)$ has been defined for $0\le
t<S_k$ and let
 \beqlb\label{2.3}
\xi = y(S_k-) + \int_{\{S_k\}}\int_{U_1} g_1(y(S_k-),u) N_1(ds,du).
 \eeqlb
By the assumption there is also a strong solution $\{x_k(t)\}$ to
 \beqlb\label{2.4}
x(t)
 \ar=\ar
\xi + \int_0^t \sigma(x(s))dB(S_k+s) + \int_0^t \int_{U_0}
g_0(x(s-),u) \tilde{N}_0(S_k+ds,du) \nnm \\
 \ar \ar
+ \int_0^t b(x(s))ds + \int_0^t\int_{U_2} g_1(x(s-),u)
N_1(S_k+ds,du).
 \eeqlb
For $S_k\le t<S_{k+1}$ we set $y(t)=x_k(t-S_k)$. By (\ref{2.2}) and
(\ref{2.4}) it is not hard to show that $\{y(t)\}$ is a strong
solution to (\ref{2.1}). On the other hand, if $\{y(t)\}$ is a
solution to (\ref{2.1}), it satisfies (\ref{2.2}) for $0\le t<S_1$
and satisfies (\ref{2.4}) for $S_k\le t<S_{k+1}$ with $\xi$ given by
(\ref{2.3}). Then the pathwise uniqueness for (\ref{2.1}) follows
from that for (\ref{2.2}) and (\ref{2.4}). \qed

Thanks to Proposition~\ref{p2.2} we may focus on the existence and
uniqueness of strong solutions of equation (\ref{2.2}). The
following proposition gives an estimate for the first moment of the
solution of (\ref{2.1}).

\bproposition\label{p2.3} Suppose that conditions (2.a,b) hold. Let
$\{x(t)\}$ be a non-negative solution of (\ref{2.2}) and let $\tau_m
= \inf\{t\ge0: x(t)\ge m\}$ for $m\ge 1$. Then $\tau_m\to \infty$
almost surely as $m\to \infty$ and
 \beqlb\label{2.5}
\mbf{E}[1+x(t)] \le \mbf{E}[1+x(0)] \exp\{Kt\}, \qquad t\ge 0.
 \eeqlb
\eproposition

\proof It is easy to show that
 \beqnn
\ar \ar\mbf{E}\bigg[\bigg|\int_0^{t\land\tau_m}\int_{U_0}
g_0(x(s-),u)1_{\{|g_0(x(s-),u)|> 1\}}\tilde{N}_0(ds,du)\bigg|\bigg] \\
 \ar \ar\qquad
\le \mbf{E}\bigg[\int_0^{t\land\tau_m}ds\int_{U_0} |g_0(x(s-),u)|
1_{\{|g_0(x(s-),u)|> 1\}} \mu_0(du)\bigg].
 \eeqnn
Observe also that
 \beqnn
\ar\ar\mbf{E}\bigg[\bigg|\int_0^{t\land\tau_m} \sigma(x(s-))
dB(s)\bigg|^2\bigg]
 =
\mbf{E}\bigg[\int_0^{t\land\tau_m} \sigma(x(s-))^2ds\bigg]
 \eeqnn
and
 \beqnn
\ar\ar\mbf{E}\bigg[\bigg|\int_0^{t\land\tau_m}\int_{U_0}
g_0(x(s-),u) 1_{\{|g_0(x(s-),u)|\le 1\}}
\tilde{N}_0(ds,du)\bigg|^2\bigg] \\
 \ar\ar\qquad
= \mbf{E}\bigg[\int_0^{t\land\tau_m}ds\int_{U_0} g_0(x(s-),u)^2
1_{\{|g_0(x(s-),u)|\le 1\}} \mu_0(du)\bigg].
 \eeqnn
Since $x(s-)\le m$ for all $0< s\le \tau_m$, the above expectations
are finite by (2.b). Consequently,
 \beqnn
t\mapsto \int_0^{t\land\tau_m} \sigma(x(s-)) dB(s) + \int_0^{t\land
\tau_m} \int_{U_0}g_0(x(s-),u) \tilde{N}_0(ds,du)
 \eeqnn
is a martingale. {From} (\ref{2.2}) and (2.a) we get
 \beqnn
\mbf{E}[1+x(t\land\tau_m)]
 \ar=\ar
\mbf{E}[1+x(0)] + \mbf{E}\bigg[\int_0^{t\land\tau_m} |b(x(s-))|
ds\bigg] \\
 \ar \ar
+\, \mbf{E}\bigg[\int_0^{t\land\tau_m}ds\int_{U_2} |g_1(x(s-),u)|
\mu_1(du)\bigg] \\
 \ar\le\ar
\mbf{E}[1+x(0)] + K\mbf{E}\bigg[\int_0^{t\land\tau_m} (1+x(s-))
ds\bigg].
 \eeqnn
Thus $t\mapsto \mbf{E}[1+x(t\land\tau_m)]$ is a locally bounded
function. Moreover, since $x(s-) \neq x(s)$ for at most countably
many $s\ge 0$, it follows that
 \beqnn
\mbf{E}[1+x(t\land\tau_m)]
 \le
\mbf{E}[1+x(0)] + K\int_0^t \mbf{E}[1+x(s\land\tau_m)] ds.
 \eeqnn
By Gronwall's lemma,
 \beqlb\label{2.6}
\mbf{E}[1+x(t\land \tau_m)]
 \le
\mbf{E}[1+x(0)] \exp\{Kt\}, \qquad t\ge 0.
 \eeqlb
By the right continuity of $\{x(t)\}$ we have $x(\tau_m)\ge m$, so
the above inequality implies
 \beqnn
(1+m)\mbf{P}\{\tau_m\le t\}
 \le
\mbf{E}[1+x(0)] \exp\{Kt\}.
 \eeqnn
Then $\tau_m\to \infty$ almost surely as $m\to \infty$, and
(\ref{2.5}) follows from (\ref{2.6}) by an application of Fatou's
lemma. \qed

\bproposition\label{p2.4} Suppose that conditions (2.a,b) hold and
for each $m\ge1$ there is a unique strong solution to
 \beqlb\label{2.7}
x(t)
 \ar=\ar
x(0) + \int_0^t \sigma(x(s)\land m)dB(s) + \int_0^t
b_m(x(s)\land m)ds \nnm \\
 \ar \ar
+ \int_0^t\int_{U_0} [g_0(x(s-)\land m,u)\land m]
\tilde{N}_0(ds,du) \nnm \\
 \ar \ar
+ \int_0^t\int_{U_2} [g_1(x(s-)\land m,u)\land m] N_1(ds,du),
 \eeqlb
where
 \beqnn
b_m(x) = b(x) - \int_{U_0} [g_0(x,u) - g_0(x,u)\land m] \mu_0(du).
 \eeqnn
Then there is a unique strong solution to (\ref{2.2}). \eproposition

\proof For any $m\ge 1$ let $\{x_m(t)\}$ denote the unique strong
solution to (\ref{2.7}) and let $\tau_m = \inf\{t\ge 0: x_m(t)\ge
m\}$. Since $0\le x_m(t)< m$ for $0\le t< \tau_m$, the trajectory
$t\mapsto x_m(t)$ has no jumps larger than $m$ on the time interval
$[0,\tau_m)$. Then for $0\le t< \tau_m$ we have
 \beqnn
x_m(t)
 \ar=\ar
x(0) + \int_0^t \sigma(x_m(s))dB(s) + \int_0^t\int_{U_0}
g_0(x_m(s-),u) \tilde{N}_0(ds,du) \nnm \\
 \ar \ar
+ \int_0^tds\int_{U_2}[g_1(x_m(s-),u) - g_1(x_m(s-),u)\land m]
\mu_1(du) \nnm \\
 \ar \ar
+ \int_0^t b_m(x_m(s))ds + \int_0^t\int_{U_2} g_1(x_m(s-),u)
N_1(ds,du) \nnm \\
 \ar=\ar
x(0) + \int_0^t \sigma(x_m(s))dB(s) + \int_0^t\int_{U_0}
g_0(x_m(s-),u) \tilde{N}_0(ds,du) \nnm \\
 \ar \ar
+ \int_0^t b(x_m(s))ds + \int_0^t\int_{U_2} g_1(x_m(s-),u)
N_1(ds,du).
 \eeqnn
In other words, $\{x_m(t)\}$ satisfies (\ref{2.2}) for $0\le t<
\tau_m$. For $n\ge m\ge 1$ let $\{y(t): t\ge 0\}$ be the unique
solution to
 \beqnn
y(t)
 \ar=\ar
\xi + \int_0^t \sigma(y(s)\land n)dB(\tau_m+s) + \int_0^t
b_n(y(s)\land n)ds \nnm \\
 \ar \ar
+ \int_0^t \int_{U_0} [g_0(y(s-)\land n,u)\land n]
\tilde{N}_0(\tau_m+ds,du) \nnm \\
 \ar \ar
+ \int_0^t\int_{U_2} [g_1(y(s-)\land n,u)\land n] N_1(\tau_m+ds,du),
 \eeqnn
where
 \beqnn
\xi \ar=\ar x_m(\tau_m-) + \int_{\{\tau_m\}}\int_{U_0}
[g_0(x_m(\tau_m-),u) \land n] N_0(ds,du) \nnm \\
 \ar \ar
+ \int_{\{\tau_m\}}\int_{U_2} [g_1(x_m(\tau_m-),u) \land n]
N_1(ds,du).
 \eeqnn
We define the process $\{x_n^\prime(t)\}$ by
 \beqnn
x_n^\prime(t)
 =
\bigg\{\begin{array}{ll}
x_m(t) \quad \ar\mbox{if $0\le t< \tau_m$,} \\
y(t-\tau_m) \ar\mbox{if $t\ge \tau_m$.}
\end{array}
 \eeqnn
It is not hard to see that $\{x_n^\prime(t)\}$ is a solution to
(\ref{2.7}) with the $m$ replaced by $n$. By the strong uniqueness
we get $x_n^\prime(t) = x_n(t)$ for all $t\ge0$. In particular, we
infer $x_n(t) = x_m(t)< m$ for $0\le t< \tau_m$. Consequently, the
sequence $\{\tau_m\}$ is non-decreasing. On the other hand, as in
the proof of Proposition~\ref{p2.3} we have
 \beqnn
\mbf{E}[1+x_m(t\land \tau_m)]
 \le
\mbf{E}[1+x(0)] \exp\{Kt\}, \qquad t\ge 0,
 \eeqnn
and hence $\tau_m\to \infty$ as $m\to \infty$. Let $\{x(t)\}$ be the
process such that $x(t) = x_m(t)$ for all $0\le t< \tau_m$ and $m\ge
1$. It is easily seen that $\{x(t)\}$ is a strong solution of
(\ref{2.2}). The uniqueness of solution follows by a similar
localization argument. \qed

Now we give a simple criterion for the existence and uniqueness of
the strong solution to (\ref{2.1}). We shall assume the following
local Lipschitz condition:
 \begin{itemize}

\item[{\rm(2.c)}]
For each integer $m\ge1$ there is a constant $K_m\ge 0$ so that
 \beqnn
|\sigma(x)-\sigma(y)| + |b(x)-b(y)| + \int_{U_2} |g_1(x,u)-g_1(y,u)|
\mu_1(du) \le K_m|x-y|
 \eeqnn
and
 \beqlb\label{2.8}
|g_0(x,u) - g_0(y,u))|\le K_m|x-y|f_m(u)
 \eeqlb
for $0\le x,y\le m$ and $u\in U_0$, where $u\mapsto f_m(u)$ is a
non-negative function on $U_0$ satisfying
 \beqnn
\int_{U_0} [f_m(u)\land f_m(u)^2] \mu_0(du)< \infty.
 \eeqnn

 \end{itemize}

\btheorem\label{t2.5} Suppose that conditions (2.a,b,c) are
satisfied. Then there is a unique non-negative strong solution to
(\ref{2.1}). \etheorem

\proof By Proposition~\ref{p2.2} we only need to show there is a
unique strong solution to (\ref{2.2}). For $m\ge 1$ let $V_m =
\{u\in U_0: f_m(u)\le 1/K_m\}$. Then (2.c) implies
$\mu_0(U_0\setminus V_m)< \infty$. Let us consider the equation
 \beqlb\label{2.9}
x(t)
 \ar=\ar
x(0) + \int_0^t \sigma(x(s)\land m)dB(s) + \int_0^t
b(x(s)\land m)ds \nnm \\
 \ar \ar
+ \int_0^t\int_{U_0} [g_0(x(s-)\land m,u)\land m]
\tilde{N}_0(ds,du) \nnm \\
 \ar \ar
- \int_0^tds\int_{U_0\setminus V_m} [g_0(x(s-)\land m,u) -
g_0(x(s-)\land m,u)\land m] \mu_0(du) \nnm \\
 \ar \ar
+ \int_0^t\int_{U_2} [g_1(x(s-)\land m,u)\land m] N_1(ds,du),
 \eeqlb
which can be rewritten as a jump-type equation with compensated
Poisson integral over $V_m\cup U_2$, non-compensated Poisson
integral over $U_0\setminus V_m$, and drift coefficient
 \beqnn
x \mapsto b(x\land m) - \int_{U_0\setminus V_m} g_0(x\land m,u)
\mu_0(du) + \int_{U_2} [g_1(x\land m,u)\land m] \mu_1(du).
 \eeqnn
By a classical result it follows that (\ref{2.9}) has a unique
strong solution; see, e.g., Ikeda and Watanabe (1989, pp.244-245).
Taking $y=0$ in (\ref{2.8}) gives $g_0(x,u)\le K_mxf_m(u)$ for $0\le
x\le m$. Then $g_0(x\land m,u)> m$ implies $f_m(u)> 1/K_m$, and
hence
 \beqnn
b_m(x\land m) = b(x\land m) - \int_{U_0\setminus V_m} [g_0(x\land
m,u) - g_0(x\land m,u)\land m] \mu_0(du).
 \eeqnn
Consequently, (\ref{2.9}) is equivalent to (\ref{2.7}). By
Proposition~\ref{p2.4} there is a unique strong solution to
(\ref{2.2}). \qed

%%%%%%%%%%%%%%%%%%%%%%%%%%%%%%%%%%%%%%%%%%%%%%

\section{Pathwise uniqueness: Non-Lipschitz conditions}

\setcounter{equation}{0}

In this section, we prove some results on the pathwise uniqueness of
solutions to (\ref{2.1}) under non-Lipschitz conditions. Suppose
that $(\sigma, b, g_0, g_1, \mu_0,\mu_1)$ are given as in the last
section. Given a function $f$ defined on a subset of $\mbb{R}$, we
note
 \beqlb\label{3.1}
\Delta_zf(x) = f(x+z) - f(x)
 \quad\mbox{and}\quad
D_zf(x) = \Delta_zf(x) - f^\prime(x)z
 \eeqlb
if the right hand sides are meaningful. Let us consider the
following condition:
 \begin{itemize}

\item[{\rm(3.a)}]
There is a continuous non-decreasing function $x\mapsto b_2(x)$ on
$\mbb{R}_+$ so that $b(x) = b_1(x) - b_2(x)$, and for each integer
$m\ge1$ there is a non-decreasing and concave function $z\mapsto
r_m(z)$ on $\mbb{R}_+$ such that $\int_{0+} r_m(z)^{-1} dz$ $=$
$\infty$ and
 \beqnn
|b_1(x)-b_1(y)| + \int_{U_2} |g_1(x,u) - g_1(y,u)|\mu_1(du)
 \le
r_m(|x-y|)
 \eeqnn
for all $0\le x,y\le m$.

 \end{itemize}

\bproposition\label{p3.1} Suppose that conditions (2.a,b) and (3.a)
hold. Then the pathwise uniqueness of solution holds for (\ref{2.2})
if for each integer $m\ge1$ there is a sequence of non-negative and
twice continuously differentiable functions $\{\phi_k\}$ with the
following properties:
 \begin{itemize}

\item[{(i)}]
$\phi_k(x) \to |x|$ non-decreasingly as $k\to \infty$;

\item[{(ii)}]
$0\le \phi_k^\prime(x) \le 1$ for $x\ge 0$ and $-1\le
\phi_k^\prime(x) \le 0$ for $x\le 0$;

\item[{(iii)}]
$\phi_k^{\prime\prime}(x) \ge 0$ for $x\in \mbb{R}$ and for each
$m\ge 1$,
 \beqnn
\phi_k^{\prime\prime}(x-y)[\sigma(x) - \sigma(y)]^2 \to 0 \qquad
(k\to \infty)
 \eeqnn
uniformly in $0\le x,y\le m$;

\item[{(iv)}]
for each $m\ge 1$,
 \beqnn
\int_{U_0} D_{l(x,y;u)}\phi_k(x-y) \mu_0(du) \to 0 \qquad (k\to
\infty)
 \eeqnn
uniformly in $0\le x,y\le m$, where $l(x,y;u) = g_0(x,u) -
g_0(y,u)$.

 \end{itemize}
\eproposition

\proof Suppose that $\{x_1(t)\}$ and $\{x_2(t)\}$ be two solutions
of (\ref{2.2}) with deterministic initial values.
Proposition~\ref{p2.3} implies that $t \mapsto \mbf{E}[x_1(t)] +
\mbf{E}[x_2(t)]$ is locally bounded. Let $\zeta(t) = x_1(t) -
x_2(t)$ for $t\ge 0$. {From} (\ref{2.2}) we have
 \beqlb\label{3.2}
\zeta(t)
 \ar=\ar
\zeta(0) + \int_0^t[\sigma(x_1(s-))-\sigma(x_2(s-))]dB(s) \nnm \\
 \ar \ar
+ \int_0^t\int_{U_0} [g_0(x_1(s-),u) - g_0(x_2(s-),u)]
\tilde{N}_0(ds,du) \nnm \\
 \ar \ar
+ \int_0^t [b(x_1(s-)) - b(x_2(s-))] ds \nnm \\
 \ar \ar
+ \int_0^t\int_{U_2} [g_1(x_1(s-),u) - g_1(x_2(s-),u)] N_1(s,du).
 \eeqlb
Let $\tau_m = \inf\{t\ge 0: x_1(t)\ge m$ or $x_2(t)\ge m\}$ for
$m\ge 1$. By (\ref{3.2}) and It\^o's formula,
 \beqlb\label{3.3}
\phi_k(\zeta(t\land\tau_m))
 \ar=\ar
\phi_k(\zeta(0)) + \int_0^{t\land\tau_m} \phi_k^\prime(\zeta(s-))
[b(x_1(s-)) - b(x_2(s-))] ds \nnm \\
 \ar \ar
+\,\frac{1}{2}\int_0^{t\land\tau_m}\phi_k^{\prime\prime}(\zeta(s-))
[\sigma(x_1(s-)) - \sigma(x_2(s-))]^2ds \nnm \\
 \ar \ar
+ \int_0^{t\land\tau_m}\int_{U_2} \Delta_{l_1(s-,u)}
\phi_k(\zeta(s-)) N_1(ds,du) \nnm \\
 \ar \ar
+ \int_0^{t\land\tau_m}\int_{U_0} D_{l_0(s-,u)}\phi_k(\zeta(s-))
N_0(ds,du) + \mbox{mart.} \nnm \\
 \ar=\ar
\phi_k(\zeta(0)) + \int_0^{t\land\tau_m} \phi_k^\prime(\zeta(s-))
[b(x_1(s-)) - b(x_2(s-))] ds \nnm \\
 \ar \ar
+\,\frac{1}{2}\int_0^{t\land\tau_m} \phi_k^{\prime\prime}(\zeta(s-))
[\sigma(x_1(s-))-\sigma(x_2(s-))]^2 ds \nnm \\
 \ar \ar
+ \int_0^{t\land\tau_m} ds \int_{U_2} \Delta_{l_1(s-,u)}
\phi_k(\zeta(s-)) \mu_1(du) \nnm \\
 \ar \ar
+ \int_0^{t\land\tau_m} ds \int_{U_0} D_{l_0(s-,u)}
\phi_k(\zeta(s-)) \mu_0(du) + \mbox{mart.},
 \eeqlb
where
 \beqnn
l_i(s,u) = g_i(x_1(s),u) - g_i(x_2(s),u), \qquad i=1,2.
 \eeqnn
Since $b(x) = b_1(x) - b_2(x)$ and $x\mapsto b_2(x)$ is
non-decreasing, by property (ii) we have
 \beqnn
\phi_k^\prime(\zeta(s-))\big[b(x_1(s-)) - b(x_2(s-))\big]
 \ar\le\ar
\phi_k^\prime(\zeta(s-))\big[b_1(x_1(s-)) - b_1(x_2(s-))\big] \\
 \ar\le\ar
|b_1(x_1(s-))-b_1(x_2(s-))|.
 \eeqnn
Observe also that
 \beqnn
\int_{U_2}\Delta_{l_1(s-,u)}\phi_k(\zeta(s-))\mu_1(du) \le
\int_{U_2} |g_1(x_1(s-),u) - g_1(x_2(s-),u)| \mu_1(du).
 \eeqnn
By condition (3.a), for any $s\le\tau_m$ the summation of the right
hand sides of the above two inequalities is no larger than
$r_m(|\zeta(s-)|)$. By properties (iii) and (iv) we have
 \beqnn
\phi_k^{\prime\prime}(\zeta(s-))[\sigma(x_1(s-)) -
\sigma(x_2(s-))]^2 \to 0
 \eeqnn
and
 \beqnn
\int_{U_0} D_{l_0(s-,u)}\phi_k(\zeta(s-)) \mu_0(du) \to 0
 \eeqnn
uniformly on the event $\{s\le\tau_m\}$. Then we can take the
expectations in (\ref{3.3}) and let $k\to \infty$ to get
 \beqnn
\mbf{E}[|\zeta(t\land\tau_m)|]
 \le
\mbf{E}[|\zeta(0)|] + \mbf{E}\bigg[\int_0^{t\land \tau_m}
r_m(|\zeta(s-)|)ds\bigg].
 \eeqnn
Since $\zeta(s-)< m$ for $0< s\le \tau_m$, we infer that $t\mapsto
\mbf{E}[|\zeta(t\land\tau_m)|]$ is locally bounded. Note also that
$\zeta(s-) \neq \zeta(s)$ for at most countably many $s\ge 0$. Then
the concaveness of $x\mapsto r_m(x)$ implies
 \beqnn
\mbf{E}[|\zeta(t\land\tau_m)|]
 \ar\le\ar
\mbf{E}[|\zeta(0)|] + \int_0^t \mbf{E}[r_m(|\zeta(s\land
\tau_m)|)]ds \\
 \ar\le\ar
\mbf{E}[|\zeta(0)|] + \int_0^t r_m(\mbf{E}[|\zeta(s\land
\tau_m)|])ds.
 \eeqnn
If $x_1(0) = x_2(0)$, we can use a standard argument to show
$\mbf{E}[|\zeta(t\land \tau_m)|] =0$ for all $t\ge0$; see e.g.\
Ikeda and Watanabe (1989, p.184). Since $\tau_m\to \infty$ as
$m\to\infty$ by Proposition~\ref{p2.3}, the right continuity of
$t\mapsto \zeta(t)$ implies $\mbf{P}\{\zeta(t) = 0$ for all $t\ge0\}
=1$. \qed

We remark that the proof of Proposition~\ref{p3.1} given above uses
essentially the monotonicity of $x\mapsto b_2(x)$. A similar
condition has been used in Rozovsky (1980) for continuous type
equations. The reader may refer Cherny and Engelbert (2005) for a
thorough treatments of continuous type stochastic differential
equations with singular coefficients. For the next theorem we need
the following condition:

 \begin{itemize}

\item[{\rm(3.b)}]
For every fixed $u\in U_0$ the function $x \mapsto g_0(x,u)$ is
non-decreasing, and for each integer $m\ge 1$ there is a
non-negative and non-decreasing function $z\mapsto \rho_m(z)$ on
$\mbb{R}_+$ so that $\int_{0+} \rho_m(z)^{-2}dz = \infty$ and
 \beqnn
|\sigma(x)-\sigma(y)|^2 + \int_{U_0} |l(x,y;u)|\land |l(x,y;u)|^2
\mu_0(du) \le \rho_m(|x-y|)^2
 \eeqnn
for all $0\le x,y\le m$, where $l(x,y;u) = g_0(x,u) - g_0(y,u)$.

 \end{itemize}

\btheorem\label{t3.2} Suppose that conditions (2.a,b) and (3.a,b)
are satisfied. Then the pathwise uniqueness of solution holds for
(\ref{2.1}). \etheorem

\proof By Proposition~\ref{p2.2} it suffices to prove the pathwise
uniqueness for (\ref{2.2}). For each integer $m\ge1$ we shall
construct a sequence of functions $\{\phi_k\}$ that satisfies the
properties required in Proposition~\ref{p3.1}. Let $a_0=1$ and
choose $a_k\to 0+$ decreasingly so that $\int_{a_k} ^{a_{k-1}}
\rho_m(z)^{-2}dz = k$ for $k\ge 1$. Let $z\mapsto \psi_k(z)$ be a
non-negative continuous function on $\mbb{R}$ which has support in
$(a_k, a_{k-1})$ and satisfies $\int_{a_k}^{a_{n-1}} \psi_k(z) dz
=1$ and $0\le \psi_k(z) \le 2 k^{-1}\rho_m(z)^{-2}$ for $a_k< z<
a_{k-1}$. For each $k\ge 1$ we define the non-negative and twice
continuously differentiable function
 \beqlb\label{3.4}
\phi_k(x) = \int_0^{|x|}dy \int_0^y\psi_k(z)dz, \qquad x\in\mbb{R}.
 \eeqlb
Clearly, the sequence $\{\phi_k\}$ satisfies properties (i) and (ii)
in Proposition~\ref{p3.1}. By condition (3.b) we have
 \beqnn
\phi_k^{\prime\prime}(x-y)[\sigma(x) - \sigma(y)]^2
 \le
\psi_k(|x-y|)\rho_m(|x-y|)^2 \le 2/k
 \eeqnn
for $0\le x,y\le m$. Thus $\{\phi_k\}$ also satisfies property
(iii). By Taylor's expansion,
 \beqnn
D_h\phi_k(\zeta)
 =
h^2\int_0^1\phi_k^{\prime\prime}(\zeta+th)(1-t)dt
 =
h^2\int_0^1\psi_k(|\zeta+th|)(1-t)dt.
 \eeqnn
Consequently, the monotonicity of $z\mapsto \rho_m(z)$ implies
 \beqlb\label{3.5}
D_h\phi_k(\zeta)
 \le
2k^{-1}h^2\int_0^1\rho_m(|\zeta+th|)^{-2}(1-t)dt
 \le
k^{-1}h^2\rho_m(|\zeta|)^{-2}
 \eeqlb
if $\zeta h\ge 0$. Observe also that
 \beqlb\label{3.6}
D_h\phi_k(\zeta)
 =
\Delta_h\phi_k(\zeta) - \phi_k^\prime(\zeta)h
 \le
\Delta_h\phi_k(\zeta)
 \le
|h|
 \eeqlb
if $\zeta h\ge 0$. Since $x\mapsto g_0(x,u)$ is non-decreasing, for
$0\le x,y\le m$ and $n\ge1$ we can use (\ref{3.5}) and (\ref{3.6})
to get
 \beqnn
\int_{U_0} D_{l(x,y;u)}\phi_k(x-y) \mu_0(du)
 \ar\le\ar
\frac{1}{k\rho_m(|x-y|)^2}\int_{U_0} l(x,y;u)^2
1_{\{|l(x,y;u)|\le n\}} \mu_0(du) \\
 \ar \ar
+ \int_{U_0} |l(x,y;u)|1_{\{|l(x,y;u)|> n\}} \mu_0(du) \\
 \ar\le\ar
\frac{n}{k\rho_m(|x-y|)^2}\int_{U_0} |l(x,y;u)|\land l(x,y;u)^2\mu_0(du) \\
 \ar \ar
+ \int_{U_0} g_0(m,u)1_{\{g_0(m,u)> n\}} \mu_0(du) \\
 \ar\le\ar
\frac{n}{k} + \int_{U_0} g_0(m,u)1_{\{g_0(m,u)> n\}} \mu_0(du),
 \eeqnn
where the last inequality follows by (3.b). From (2.b) we see that
$\{\phi_k\}$ satisfies property (iv) in Proposition~\ref{p3.1}. That
proves the pathwise uniqueness for (\ref{2.2}). \qed

The following conditions on the modulus of continuity are
particularly useful in applications to stochastic equations driven
by L\'evy processes:
 \begin{itemize}

\item[{\rm(3.c)}]
For every fixed $u\in U_0$ the function $x \mapsto g_0(x,u)$ is
non-decreasing, and for each integer $m\ge1$ there is a non-negative
and non-decreasing function $z\mapsto \rho_m(z)$ on $\mbb{R}_+$ so
that $\int_{0+} \rho_m(z)^{-2}dz = \infty$,
 \beqnn
|\sigma(x)-\sigma(y)| \le \rho_m(|x-y|)
 \quad\mbox{and}\quad
|g_0(x,u)-g_0(y,u)|\le \rho_m(|x-y|)f_m(u)
 \eeqnn
for all $0\le x,y\le m$ and $u\in U_0$, where $u\mapsto f_m(u)$ is a
non-negative function on $U_0$ satisfying
 \beqnn
\int_{U_0} [f_m(u)\land f_m(u)^2] \mu_0(du)<\infty.
 \eeqnn

 \end{itemize}

\btheorem\label{t3.3} Suppose that conditions (2.a,b) and (3.a,c)
are satisfied. Then the pathwise uniqueness for (\ref{2.1}) holds.
\etheorem

\proof By Proposition~\ref{p2.2} we only need to show the pathwise
uniqueness for (\ref{2.2}). The first part of this proof is
identical with that of Theorem~\ref{t3.2}. Under condition (3.c) we
have
 \beqnn
D_{l(x,y;u)}\phi_k(x-y)
 \le
k^{-1}l(x,y;u)^2\rho_m(|x-y|)^{-2}
 \le
k^{-1}f_m(u)^2
 \eeqnn
and
 \beqnn
D_{l(x,y;u)}\phi_k(x-y)
 \le
|l(x,y;u)|
 \le
\rho_m(|x-y|)f(u)
 \le
\rho_m(m)f(u)
 \eeqnn
for all $0\le x,y\le m$ and $u\in U_0$. By (3.c) for any $n\ge 1$ we
have
 \beqnn
\int_{U_0} D_{l(x,y;u)}\phi_k(x-y) \mu_0(du)
 \ar\le\ar
\frac{1}{k}\int_{U_0} f_m(u)^21_{\{f_m(u)\le n\}} \mu_0(du) \\
 \ar \ar
+\, \rho_m(m)\int_{U_0} f_m(u)1_{\{f_m(u)>n\}} \mu_0(du) \\
 \ar\le\ar
\frac{n}{k}\int_{U_0} [f_m(u)\land f_m(u)^2] \mu_0(du) \\
 \ar \ar
+\, \rho_m(m)\int_{U_0} f_m(u)1_{\{f_m(u)>n\}} \mu_0(du).
 \eeqnn
By letting $k\to \infty$ and $n\to \infty$ we see that $\{\phi_k\}$
satisfies property (iv) in Proposition~\ref{p3.1}. Then the pathwise
uniqueness for (\ref{2.2}) holds. \qed

%%%%%%%%%%%%%%%%%%%%%%%%%%%%%%%%%%%%%%%%%%%%%%

\section{Martingale problems and weak solutions}

\setcounter{equation}{0}

In this section, we study the existence of weak solutions of
(\ref{2.1}). As for continuous-type equations, this is closely
related with the corresponding martingale problems. We define the
operator $A$ from $C^2(\mbb{R}_+)$ to $C(\mbb{R}_+)$ by
 \beqlb\label{4.1}
Af(x) \ar=\ar \frac{1}{2}\sigma(x)^2f^{\prime\prime}(x) +
b(x)f^\prime(x) + \int_{U_0}D_{g_0(x,u)}f(x) \mu_0(du) \nnm \\
 \ar\ar
+ \int_{U_1} \Delta_{g_1(x,u)}f(x) \mu_1(du).
 \eeqlb
To simplify the statements we introduce the following conditions:
 \begin{itemize}

\item[(4.a)]
There is a constant $K\ge 0$ such that
 \beqnn
\ar\ar \sup_{x\ge 0} \big[b(x)^2 + \sigma(x)^2\big] + \int_{U_0}
\sup_{x\ge 0}g_0(x,u)^2 \mu_0(du) \nnm \\
 \ar\ar\qquad\qquad
+ \int_{U_0} \sup_{x\ge 0} \big[|g_1(x,u)|\vee g_1(x,u)^2\big]
\mu_1(du) \le K;
 \eeqnn

\item[(4.b)]
$x \mapsto g_0(x,u)$ is non-decreasing and continuous in measure
with respect to $\mu_0(du)$;

\item[(4.c)]
$x \mapsto g_1(x,u)$ is continuous in measure with respect to
$\mu_1(du)$.

 \end{itemize}

\bproposition\label{p4.1} If condition (4.a) holds, for any solution
$\{x(t)\}$ of (\ref{2.1}) we have
 \beqlb\label{4.2}
\mbf{E}\Big[\sup_{0\le s\le t}x(s)^2\Big]
 \le
6\mbf{E}[x(0)^2] + 24Kt + 6K^2t^2, \qquad t\ge 0.
 \eeqlb
\eproposition

\proof We first write (\ref{2.1}) into
 \beqlb\label{4.3}
x(t)
 \ar=\ar
x(0) + \int_0^t \sigma(x(s))dB(s) + \int_0^t\int_{U_0} g_0(x(s-),u)
\tilde{N}_0(ds,du) \nnm \\
 \ar \ar
+ \int_0^t b(x(s))ds + \int_0^t\int_{U_1} g_1(x(s-),u)
\tilde{N}_1(ds,du) \nnm \\
 \ar \ar
+ \int_0^tds\int_{U_1} g_1(x(s),u) \mu_1(du).
 \eeqlb
By applying Doob's inequality to the martingale terms in (\ref{4.3})
we have
 \beqnn
\mbf{E}\Big[\sup_{0\le s\le t}x(s)^2\Big]
 \ar\le\ar
6\mbf{E}[x(0)^2] + 24\mbf{E}\bigg[\int_0^t\sigma(x(s))^2ds\bigg] +
6\mbf{E}\bigg[\bigg(\int_0^t b(x(s)) ds\bigg)^2\bigg] \\
 \ar \ar
+\, 24\mbf{E}\bigg[\int_0^tds\int_{U_0}g_0(x(s),u)^2\mu_0(du)\bigg] \\
 \ar \ar
+\, 24\mbf{E}\bigg[\int_0^tds\int_{U_1}g_1(x(s),u)^2\mu_1(du)\bigg] \\
 \ar \ar
+\, 6\mbf{E}\bigg[\bigg(\int_0^tds\int_{U_1}g_1(x(s),u)
\mu_1(du)\bigg)^2\bigg] \\
 \ar\le\ar
6\mbf{E}[x(0)^2] + 24Kt + 6t\mbf{E}\bigg[\int_0^t b(x(s))^2 ds\bigg] \\
 \ar \ar
+\, 6t\mbf{E}\bigg[\int_0^t\bigg(\int_{U_1}g_1(x(s),u)
\mu_1(du)\bigg)^2ds\bigg] \\
 \ar\le\ar
6\mbf{E}[x(0)^2] + 24Kt + 6K^2t^2.
 \eeqnn
That proves (\ref{4.2}). \qed

\bproposition\label{p4.2} Suppose that condition (4.a) hold. Then a
non-negative c\`adl\`ag process $\{x(t)\}$ is a weak solution of
(\ref{2.1}) if and only if for every $f\in C^2(\mbb{R}_+)$,
 \beqlb\label{4.4}
f(x(t)) - f(x(0)) - \int_0^tAf(x(s))ds, \qquad t\ge 0
 \eeqlb
is a martingale. \eproposition

\proof Without loss of generality, we assume $x(0)$ is
deterministic. If $\{x(t)\}$ is a solution of (\ref{2.1}), by
It\^o's formula it is easy to see that (\ref{4.4}) is a bounded
martingale. Conversely, suppose that (\ref{4.4}) is a martingale for
every $f\in C^2(\mbb{R}_+)$. By a standard stopping time argument,
we have
 \beqlb\label{4.5}
x(t) = x(0) + \int_0^tb(x(s))ds + \int_0^tds\int_{U_1} g_1(x(s),u)
\mu_1(du) + M(t)
 \eeqlb
for a square-integrable martingale $\{M(t)\}$. Let $N(ds,dz)$ be the
optional random measure $[0,\infty) \times \mbb{R}$ defined on by
 \beqnn
N(ds,dz) = \sum_{s>0}1_{\{\Delta x(s)\ne 0\}}\delta_{(s,\Delta
x(s))}(ds,dz),
 \eeqnn
where $\Delta x(s) = x(s) - x(s-)$. Let $\hat N(ds,dz)$ be the
predictable compensator of $N(ds,dz)$ and let $\tilde{N}(ds,dz)$
denote the compensated random measure. By (\ref{4.5}) and
Dellacherie and Meyer (1982, p.376) we have
 \beqlb\label{4.6}
x(t) = x(0) + \int_0^tb(x(s))ds + \int_0^tds\int_{U_1} g_1(x(s),u)
\mu_1(du) + M_c(t) + M_d(t),
 \eeqlb
where $\{M_c(t)\}$ is a continuous martingale and
 \beqnn
M_d(t) = \int_0^t\int_{\mbb{R}} z \tilde N(ds,dz)
 \eeqnn
is a purely discontinuous martingale. Let $\{C(t)\}$ denote the
quadratic variation process of $\{M_c(t)\}$. By (\ref{4.6}) and
It\^o's formula we have
 \beqlb\label{4.7}
f(x(t)) \ar=\ar f(x(0)) + \int_0^t f^\prime(x(s))b(x(s))ds +
\frac{1}{2}\int_0^t f^{\prime\prime}(x(s))d C(s) \nnm \\
 \ar \ar
+ \int_0^tf^\prime(x(s))ds\int_{U_1} g_1(x(s),u) \mu_1(du) \nnm \\
 \ar \ar
+ \int_0^t\int_{\mbb{R}} D_zf(x(s-)) \hat{N}(ds,dz) + \mbox{mart.}
 \eeqlb
In view of (\ref{4.4}) and (\ref{4.7}), the uniqueness of canonical
decompositions of semi-martingales implies $dC(s) = \sigma(x(s))ds$
and
 \beqnn
\ar\ar\int_0^t\int_{\mbb{R}} F(s,u) \hat{N}(ds,dz) \\
 \ar\ar\quad
= \int_0^tds\int_{U_0} F(s,g_0(x(s-),u))\mu_0(du) +
\int_0^tds\int_{U_1} F(s,g_1(x(s-),u))\mu_1(du)
 \eeqnn
for any non-negative Borel function $F$ on $\mbb{R}_+\times
\mbb{R}$. By martingale representation theorems we have the equation
(\ref{2.1}) on an extension of the original probability space; see,
e.g., Ikeda and Watanabe (1989, p.84 and p.93). \qed

For simplicity we assume the initial variable $x(0)$ is
deterministic in the sequel of this section. To prove the existence
of a weak solution of (\ref{2.1}) let $\{\varepsilon_n\}$ be a
sequence of strictly positive numbers decreasing to zero. Let $\bar
g_i(x,u) = \sup_{0\le y\le x} |g_i(y,u)|$ for $i=0,1$. If condition
(4.b) holds,
 \beqnn
x \mapsto \int_{U_0} 1_{\{\bar g_0(n,u)> \varepsilon_n\}} g_0(x,u)
\mu_0(ds,du)
 \eeqnn
is a bounded continuous non-decreasing function on $\mbb{R}_+$. By
the result on continuous type stochastic equations, there is a weak
solution to
 \beqlb\label{4.8}
x(t)
 \ar=\ar
x(0) + \int_0^t \sigma(x(s))dB(s) + \int_0^t b(x(s))ds \nnm \\
 \ar \ar\qquad\qquad
- \int_0^tds\int_{U_0} 1_{\{\bar g_0(n,u)> \varepsilon_n\}}
g_0(x(s),u) \mu_0(du);
 \eeqlb
see, e.g., Ikeda and Watanabe (1989, p.169). By Theorem~\ref{t3.2}
the pathwise uniqueness holds for (\ref{4.8}), so the equation has a
unique strong solution. Under condition (4.a) we have
 \beqnn
\mu_i(\{u\in U_i: \bar g_i(n,u)> \varepsilon_n\})
 \le
\frac{1}{\varepsilon_n^2}\int_{U_i} \bar g_i(n,u)^2 \mu_i(du)
 \le
\frac{K}{\varepsilon_n^2}
 <
\infty.
 \eeqnn
Then by Proposition~\ref{p2.2} for every integer $n\ge 1$ there is a
unique strong solution $\{x_n(t): t\ge0\}$ to
 \beqlb\label{4.9}
x(t)
 \ar=\ar
x(0) + \int_0^t \sigma(x(s))dB(s) + \int_0^tb(x(s))ds \nnm \\
 \ar \ar
+ \int_0^t\int_{U_0} 1_{\{\bar g_0(n,u)> \varepsilon_n\}}
g_0(x(s-),u) \tilde{N}_0(ds,du) \nnm \\
 \ar \ar
+ \int_0^t\int_{U_1} 1_{\{\bar g_1(n,u)> \varepsilon_n\}}
g_1(x(s-),u) N_1(ds,du).
 \eeqlb

\blemma\label{l4.3} Under conditions (4.a,b), the sequence
$\{x_n(t): t\ge0\}$ is tight in the Skorokhod space $D([0,\infty),
\mbb{R}_+)$. \elemma

\proof Since $x(0)$ is deterministic, by Proposition~\ref{p4.1} it
is easy to see that
 \beqnn
t \mapsto C(t) := \sup_{n\ge 1}\mbf{E}\Big[\sup_{0\le s \le t}
x_n(s)^2\Big]
 \eeqnn
is locally bounded. Then for every fixed $t\ge 0$ the sequence of
random variables $x_n(t)$ is tight. Moreover, in view of
(\ref{4.3}), if $\{\tau_n\}$ is a sequence of stopping times bounded
above by $T\ge 0$, we have
 \beqnn
\mbf{E}[|x_n(\tau_n+t)-x_n(\tau_n)|^2]
 \ar=\ar
5\mbf{E}\bigg[\int_0^t\sigma(x_n(\tau_n+s))^2]ds\bigg] \nnm \\
 \ar \ar
+ 5t\mbf{E}\bigg[\int_0^tb(x_n(\tau_n+s))^2ds\bigg] \nnm \\
 \ar \ar
+\, 5\mbf{E}\bigg[\int_0^tds\int_{U_0} g_0(x_n(\tau_n+s),u)^2
\mu_0(du)\bigg] \nnm \\
 \ar \ar
+\, 5\mbf{E}\bigg[\int_0^tds\int_{U_1} g_1(x_n(\tau_n+s),u)^2
\mu_1(du)\bigg] \nnm \\
 \ar \ar
+\, 5t\mbf{E}\bigg[\int_0^t\bigg(\int_{U_1} g_1(x_n(\tau_n+s),u)
\mu_1(du)\bigg)^2ds\bigg] \nnm \\
 \ar\le\ar
5Kt(1+Kt),
 \eeqnn
where the last inequality follows by (4.a). Consequently, as $t\to
0$,
 \beqnn
\sup_{n\ge 1}\mbf{E}[|w_n(\tau_n+t)-w_n(\tau_n)|^2] \to 0.
 \eeqnn
Then $\{x_n(t): t\ge 0\}$ is tight in $D([0,\infty), \mbb{R}_+)$ by
the criterion of Aldous (1978). \qed

\btheorem\label{t4.4} Under the conditions (4.a,b,c), there exists a
non-negative weak solution to (\ref{2.1}). \etheorem

\proof Let $\{x_n(t): t\ge0\}$ be the unique non-negative strong
solution of (\ref{4.9}). By Proposition~\ref{p4.2}, for every $f\in
C^2(\mbb{R}_+)$,
 \beqlb\label{4.10}
f(x_n(t)) - f(x_n(0)) - \int_0^tA_nf(x_n(s))ds, \qquad t\ge 0
 \eeqlb
is a bounded martingale, where
 \beqnn
A_nf(x) \ar=\ar \frac{1}{2}\sigma(x)^2f^{\prime\prime}(x) +
\int_{U_0} 1_{\{\bar g_0(n,u)> \varepsilon_n\}} D_{g_0(x,u)}f(x)
\mu_0(du) \nnm \\
 \ar\ar
+\, b(x)f^\prime(x) + \int_{U_1} 1_{\{\bar g_1(n,u)>
\varepsilon_n\}} \Delta_{g_1(x,u)}f(x) \mu_1(du).
 \eeqnn
By Lemma~\ref{l4.3} there is a subsequence $\{x_{n_k}(t): t\ge0\}$
of $\{x_n(t): t\ge0\}$ that converges to some process $\{x(t):
t\ge0\}$ in distribution on $D([0,\infty), \mbb{R}_+)$. By Skorokhod
representation we may assume $\{x_{n_k}(t): t\ge0\}$ converges to
$\{x(t): t\ge0\}$ almost surely in $D([0,\infty), \mbb{R}_+)$. Let
$D(x) := \{t>0: \mbf{P}\{x(t-)\neq x(t)\}= 1\}$. Then the set
$[0,\infty) \setminus D(x)$ is at most countable; see, e.g., Ethier
and Kurtz (1986, p.131). It follows that $\lim_{k\to \infty}
x_{n_k}(t) = x(t)$ almost surely for every $t\in D(x)$; see, e.g.,
Ethier and Kurtz (1986, p.125). Using conditions (2.a,b,c) it is
elementary to show that (\ref{4.4}) is a bounded martingale. Then
the theorem follows by another application of
Proposition~\ref{p4.2}. \qed

%%%%%%%%%%%%%%%%%%%%%%%%%%%%%%%%%%%%%%%

\section{Strong solutions: Non-Lipschitz conditions}

\setcounter{equation}{0}

In this section, we give some criterion on the existence and
uniqueness of the strong solution of equation (\ref{2.1}). We also
prove two simple properties of the solution, continuous dependence
on the initial value and comparison property.

\btheorem\label{t5.1} Suppose that conditions (2.a,b) and (3.a,b)
are satisfied. Then there exists a unique non-negative strong
solution to (\ref{2.1}). \etheorem

\proof By Proposition~\ref{p2.2} we only need to (\ref{2.2}) has a
unique strong solution. Let us consider the sequence of equations
(\ref{2.7}). Since $x\mapsto g_0(x,u)$ is non-decreasing, so is
 \beqnn
x\mapsto \beta_m(x) := \int_{U_0} [g_0(x,u) - g_0(x,u)\land m]
\mu_0(du).
 \eeqnn
By conditions (3.a,b) it is easy to show that $x\mapsto g_0(x,u)$
and $x\mapsto g_1(x,u)$ are continuous in measure relative to
$\mu_0(du)$ and $\mu_1(du)$, respectively. Then one may use
dominated convergence to see $x\mapsto \beta_m(x)$ is continuous. By
applying Theorem~\ref{t4.4} we infer that (\ref{2.7}) has a weak
solution for every $m\ge 1$. The pathwise uniqueness for (\ref{2.7})
holds by Theorems~\ref{t3.2}. Thus Proposition~\ref{p2.4} implies
that (\ref{2.2}) has a unique strong solution. \qed

If (\ref{2.1}) has a unique solution, it is easy to see that the
solution is a strong Markov process with generator $A$ defined by
(\ref{4.1}). Based on the pathwise uniqueness stated in
Theorem~\ref{t3.3}, the following results can be proved similarly as
the above.

\btheorem\label{t5.2} Suppose that conditions (2.a,b) and (3.a,c)
are satisfied. Then there exists a unique non-negative strong
solution to (\ref{2.1}). \etheorem

We next present two results on the properties of the strong solution
of (\ref{2.2}). In the following two theorems we can replace
condition (3.b) by (3.c).

\btheorem\label{t5.3} Suppose that conditions (2.a,b) and (3.a,b)
hold with $U_2=U_1$ and with $r_m(x)\equiv r(x)$ independent of
$m\ge 1$. For each integer $n\ge 0$ let $\{x_n(t)\}$ be a solution
of (\ref{2.2}) and assume $\sup_{n\ge 0} \mbf{E}[x_n(0)]< \infty$.
If\, $\lim_{n\to \infty} \mbf{E}[|x_n(0) - x_0(0)|] = 0$, then
 \beqlb\label{5.1}
\lim_{n\to \infty} \mbf{E}\big[|x_n(t) - x_0(t)|\big] =0, \qquad
t\ge 0.
 \eeqlb
\etheorem

\proof Let $\zeta_n(t) = x_n(t) - x_0(t)$ and $\tau_m = \inf\{t\ge0:
\zeta_n(t)\ge m\}$. By Proposition~\ref{p2.3} it is easy to see that
$t\mapsto \mbf{E}[|\zeta_n(t)|]$ is uniformly bounded on each
bounded interval. By calculations in the proof of
Proposition~\ref{p3.1} we obtain
 \beqnn
\mbf{E}[|\zeta_n(t\land\tau_m)|]
 \ar\le\ar
\mbf{E}[|\zeta_n(0)|] + \mbf{E}\bigg[\int_0^{t\land \tau_m}
r(|\zeta_n(s-)|)ds\bigg] \\
 \ar\le\ar
\mbf{E}[|\zeta_n(0)|] + \mbf{E}\bigg[\int_0^t r(|\zeta_n(s-)|)
ds\bigg] \\
 \ar\le\ar
\mbf{E}[|\zeta_n(0)|] + \int_0^t r(\mbf{E}[|\zeta_n(s)|])ds.
 \eeqnn
Since $\tau_m\to \infty$ as $m\to \infty$, we use Fatou's lemma to
obtain
 \beqnn
\mbf{E}[|\zeta_n(t)|]
 \le
\mbf{E}[|\zeta_n(0)|] + \int_0^t r(\mbf{E}[|\zeta_n(s)|])ds.
 \eeqnn
Then another application of Fatou's lemma gives
 \beqnn
\limsup_{n\to \infty} \mbf{E}[|\zeta_n(t)|] \le \int_0^t
r\Big(\limsup_{n\to \infty} \mbf{E}[|\zeta_n(s)|]\Big)ds,
 \eeqnn
from which we obtain (\ref{5.1}). \qed

\btheorem\label{t5.4} Suppose that conditions (2.a,b) and (3.a,b)
hold. In addition, we assume the function $x\mapsto g_1(x,u)$ is
non-decreasing for every $u\in U_1$. If $\{x_1(t)\}$ and
$\{x_2(t)\}$ are solutions of (\ref{2.2}) satisfying
$\mbf{P}\{x_1(0)\le x_2(0)\} = 1$, then we have $\mbf{P} \{x_1(t)
\le x_2(t)$ for all $t\ge 0\} = 1$. \etheorem

\proof The following arguments are similar to those in the proofs of
Proposition~\ref{p3.1} and Theorem~\ref{t3.2}. Let $\zeta(t) =
x_1(t) - x_2(t)$ for $t\ge0$. Instead of (\ref{3.4}), for each
$k\ge1$ we now define
 \beqlb\label{5.2}
\phi_k(x) = \int_0^{x} dy \int_0^y\psi_k(z)dz, \qquad x\in \mbb{R}.
 \eeqlb
Then $\phi_k(x) \to x^+ := 0\vee x$ increasingly as $k\to \infty$.
Since $\phi_k(x) =0$ for $x\le 0$ and both $x\mapsto g_0(x,u)$ and
$x\mapsto g_1(x,u)$ are non-decreasing, we have
 \beqnn
\phi_k(\zeta(t\land\tau_m))
 \ar=\ar
\int_0^{t\land\tau_m} \phi_k^\prime(\zeta(s-)) [b(x_1(s-))
- b(x_2(s-))] 1_{\{\zeta(s-)>0\}}ds \nnm \\
 \ar \ar
+\,\frac{1}{2}\int_0^{t\land\tau_m} \phi_k^{\prime\prime}(\zeta(s-))
[\sigma(x_1(s-))-\sigma(x_2(s-))]^21_{\{\zeta(s-)>0\}} ds \nnm \\
 \ar \ar
+ \int_0^{t\land\tau_m}\int_{U_2} \Delta_{l_1(s-,u)}
\phi_k(\zeta(s-))1_{\{\zeta(s-)>0\}}N_1(ds,du) \nnm \\
 \ar \ar
+ \int_0^{t\land\tau_m}\int_{U_0} D_{l_0(s-,u)}\phi_k(\zeta(s-))
1_{\{\zeta(s-)>0\}} \tilde{N}_0(ds,du) + \mbox{mart.} \\
 \ar=\ar
\int_0^{t\land\tau_m} \phi_k^\prime(\zeta(s-)) [b(x_1(s-))
- b(x_2(s-))] 1_{\{\zeta(s-)>0\}}ds \nnm \\
 \ar \ar
+\,\frac{1}{2}\int_0^{t\land\tau_m} \phi_k^{\prime\prime}(\zeta(s-))
[\sigma(x_1(s-))-\sigma(x_2(s-))]^21_{\{\zeta(s-)>0\}} ds \nnm \\
 \ar \ar
+ \int_0^{t\land\tau_m} ds \int_{U_2} \Delta_{l_1(s-,u)}
\phi_k(\zeta(s-))1_{\{\zeta(s-)>0\}}\mu_1(du) \nnm \\
 \ar \ar
+ \int_0^{t\land\tau_m} ds \int_{U_0} D_{l_0(s-,u)}
\phi_k(\zeta(s-)) 1_{\{\zeta(s-)>0\}} \mu_0(du) + \mbox{mart.}
 \eeqnn
(If $x\mapsto g_1(x,u)$ were allowed to decrease, we could not
insert $1_{\{\zeta(s-)>0\}}$ into the integral with respect to
$N_1(ds,du)$.) {From} the above equation and the estimates in the
proof of Theorem~\ref{t3.2} we obtain
 \beqnn
\mbf{E}[\zeta(t\land\tau_m)^+]
 \le
\int_0^t \mbf{E}[r(\zeta(s\land\tau_m)^+)]ds
 \le
\int_0^t r(\mbf{E}[\zeta(s\land\tau_m)^+])ds.
 \eeqnn
Then $\mbf{E}[\zeta(t\land\tau_m)^+] =0$ for all $t\ge0$, giving the
desired comparison result. \qed

%%%%%%%%%%%%%%%%%%%%%%%%%%%%%%%%%%%%%%%

\section{Applications}

\setcounter{equation}{0}

In this last section we illustrate some applications of our main
results (Theorems~\ref{t2.5}, \ref{t5.1} and~\ref{t5.2}) to
construct non-negative Markov processes. Those results give great
flexibilities in the constructions. We start with a class of
non-negative processes with non-negative jumps. Suppose that
 \begin{itemize}

\item
$x\mapsto \sigma(x)$ and $x\mapsto b(x)$ are continuous functions on
$\mbb{R}_+$ satisfying $\sigma(0) =0$ and $b(0) \ge 0$;

\item
$(x,u) \mapsto h_0(x,z)$ is a non-negative Borel function on
$\mbb{R}_+ \times (0,\infty)$ so that $h_0(0,z)=0$;

\item
$(x,z) \mapsto h_1(x,z)$ is a non-negative Borel function on
$\mbb{R}_+ \times (0,\infty)$.

 \end{itemize}
Let $\mu_0(dz)$ and $\mu_1(dz)$ be $\sigma$-finite measures on
$(0,\infty)$ such that
 \beqnn
\int_0^\infty h_0(x,z)(z\land z^2) \mu_0(dz) + \int_0^\infty
h_1(x,z)(1\land z) \mu_1(dz)< \infty, \quad x\ge 0.
 \eeqnn
In view of the characterization of Feller semigroups given by
Courr\'ege (1966), it is natural consider the generator $L$ defined
by
 \beqlb\label{6.1}
Lf(x)
 \ar=\ar
\frac{1}{2}\sigma(x)^2f^{\prime\prime}(x) + \int_0^\infty [f(x+z) -
f(x) - z
f^\prime(x)] h_0(x,z)\mu_0(dz) \nnm \\
 \ar \ar
+\, b(x)f^{\prime}(x) + \int_0^\infty [f(x+z) - f(x)]
h_1(x,z)\mu_1(dz).
 \eeqlb
Let $\{B(t)\}$ be a standard $(\mcr{F}_t)$-Brownian motion, and let
$\{p_0(t)\}$ and $\{p_1(t)\}$ be $(\mcr{F}_t)$-Poisson point
processes on $(0,\infty)^2$ with characteristic measures
$\mu_0(dz)du$ and $\mu_1(dz)du$, respectively. We assume $\{B(t)\}$,
$\{p_0(t)\}$ and $\{p_1(t)\}$ are independent of each other. Let
$N_0(ds,dz,du)$ and $N_1(ds,dz,du)$ be the Poisson random measures
associated with $\{p_0(t)\}$ and $\{p_1(t)\}$, respectively. Let us
consider the stochastic equation
 \beqlb\label{6.2}
x(t)
 \ar=\ar
x(0) + \int_0^t \sigma(x(s-))dB(s) + \int_0^t\int_0^\infty
\int_0^{h_0(x(s-),z)} z \tilde{N}_0(ds,dz,du) \nnm \\
 \ar \ar
+ \int_0^t b(x(s-))ds + \int_0^t\int_0^\infty \int_0^{h_1(x(s-),z)}
z N_1(ds,dz,du).
 \eeqlb
For the ingredients of this equation we may rephrase (2.a,b) and
(3.a,b) into the following conditions:
 \begin{itemize}

\item[{\rm(6.a)}]
There is a constant $K\ge0$ so that
 \beqnn
|b(x)| + \int_0^\infty \sup_{0\le y\le x}h_1(y,z) z \mu_1(dz) \le
K(1+x), \quad x\ge 0;
 \eeqnn

\item[{\rm(6.b)}]
There is a non-negative and non-decreasing function $x\mapsto L(x)$
on $\mbb{R}_+$ such that
 \beqnn
\sigma(x)^2 + \int_0^\infty h_0(x,z)(z\land z^2) \mu_0(dz) \le L(x),
\quad x\ge 0;
 \eeqnn

\item[{\rm(6.c)}]
There is a continuous non-decreasing function $x\mapsto b_2(x)$ on
$\mbb{R}_+$ so that $b(x) = b_1(x) - b_2(x)$, and for each integer
$m\ge 1$ there is a non-decreasing concave function $u\mapsto
r_m(z)$ on $\mbb{R}_+$ such that $\int_{0+} r_m(z)^{-1}dz = \infty$
and
 \beqnn
|b_1(x)-b_1(y)| + \int_0^\infty |h_1(x,z)-h_1(y,z)| z\mu_1(dz) \le
r_m(|x-y|)
 \eeqnn
for all $0\le x,y \le m$;

\item[{\rm(6.d)}]
For fixed $z>0$ the function $x \mapsto h_0(x,z)$ is non-decreasing,
and for each integer $m\ge1$ there is a non-decreasing function
$z\mapsto \rho(z)$ on $\mbb{R}_+$ such that $\int_{0+}
\rho_m(z)^{-2}dz = \infty$ and
 \beqnn
|\sigma(x)-\sigma(y)|^2 + \int_0^\infty |h_0(x,z)-h_0(y,z)| (z\land
z^2) \mu_0(dz) \le \rho_m(|x-y|)^2
 \eeqnn
for all $0\le x,y \le m$.

 \end{itemize}

\btheorem\label{t6.1} If conditions (6,a,b,c,d) are satisfied, there
exists a unique non-negative strong solution to (\ref{6.2}) and the
solution is a strong Markov process with generator given by
(\ref{6.1}). \etheorem

\proof By applying Theorem~\ref{t5.1} with $U_0= U_1= (0,\infty)^2$
it is simple to see that (\ref{6.2}) has a unique strong solution
$\{x(t)\}$ to. The uniqueness implies the strong Markov property of
$\{x(t)\}$. By It\^o's formula it is easy to show that $\{x(t)\}$
has generator given by (\ref{6.1}). \qed

By an easy application of Theorem~\ref{t6.1} we get the following
result, which improves slightly a theorem of Dawson and Li (2006).

\bcorollary\label{c6.2} There exists a unique non-negative strong
solution to (\ref{1.7}) and the solution is CBI-process with
generator given by (\ref{1.5}). \ecorollary

Next we consider stochastic equations driven by one-sided L\`evy
processes. Let $\mu_0(dz)$ and $\mu_1(dz)$ be $\sigma$-finite
measures on $(0,\infty)$ and let $\nu_0(dz)$ and $\nu_1(dz)$ be
$\sigma$-finite measures on $(0,1]$. We assume that
 \beqnn
\int_0^\infty (z\land z^2) \mu_0(dz) + \int_0^\infty (1\land z)
\mu_1(dz) + \int_0^1 z^2 \nu_0(dz) + \int_0^1 z \nu_1(dz)< \infty.
 \eeqnn
Let $\{B(t)\}$ be a standard $(\mcr{F}_t)$-Brownian motion. Let
$\{z_0(t)\}$ and $\{z_1(t)\}$ be $(\mcr{F}_t)$-L\'evy processes with
exponents
 \beqnn
u\mapsto \int_0^\infty \big(e^{iuz} - 1 - iuz\big) \mu_0(dz)
 \quad\mbox{and}\quad
u\mapsto \int_0^\infty \big(e^{iuz} - 1\big) \mu_1(dz),
 \eeqnn
respectively.  Let $\{y_0(t)\}$ and $\{y_1(t)\}$ be
$(\mcr{F}_t)$-L\'evy processes with exponents
 \beqnn
u\mapsto \int_0^1 \big(e^{iuz} - 1 - iuz\big) \nu_0(dz)
 \quad\mbox{and}\quad
u\mapsto \int_0^1 \big(e^{iuz} - 1\big) \nu_1(dz),
 \eeqnn
respectively. Therefore $\{z_0(t)\}$ and $\{y_0(t)\}$ are centered
and $\{z_1(t)\}$ and $\{y_1(t)\}$ are non-decreasing. Suppose that
those processes are independent of each other. In addition, suppose
that
 \begin{itemize}

\item
$x\mapsto \sigma(x)$ and $x\mapsto b(x)$ are continuous functions on
$\mbb{R}_+$ satisfying $\sigma(0) =0$ and $b(0) \ge 0$;

\item
$x\mapsto \phi_0(x)$ and $x\mapsto \phi_1(x)$ are a continuous
non-negative functions on $\mbb{R}_+$ so that $\phi_0(0) =0$.

 \end{itemize}

\btheorem\label{t6.3} If the four functions introduced above are
Lipschitz, then there is a unique non-negative strong solution to
 \beqlb\label{6.3}
dx(t)
 \ar=\ar
\sigma(x(t))dB(t) + \phi_0(x(t-))dz_0(t) + b(x(t))dt  \nnm \\
 \ar \ar
+\, \phi_1(x(t-))dz_1(t) - x(t-)dy_0(t) - x(t-)dy_1(t).
 \eeqlb
\etheorem

\proof To save the notation, we extend $\mu_0$ and $\mu_1$ to
$\sigma$-finite measures on $[-1,0)\cup(0,\infty)$ by setting
$\mu_0([-x,0)) = \nu_0((0,x])$ and $\mu_1([-x,0)) = \nu_1((0,x])$
for $0<x\le 1$. By the general result on L\'evy-It\^o
decompositions, we have
 \beqnn
z_0(t) =\ \int_0^t\int_0^\infty z \tilde{N}_0(ds,dz),
 \qquad
z_1(t) = \int_0^t\int_0^\infty z N_1(ds,dz)
 \eeqnn
and
 \beqnn
y_0(t) = -\int_0^t\int_{[-1,0)} z \tilde{N}_0(ds,dz),
 \quad
y_1(t) = -\int_0^t\int_{[-1,0)} z N_1(ds,dz),
 \eeqnn
where $N_0(ds,dz)$ and $N_1(ds,dz)$ are independent Poisson random
measures with intensities $ds\mu_0(dz)$ and $ds\mu_1(dz)$,
respectively. By Theorem~\ref{t2.5} there is a unique strong
solution to
 \beqnn
x(t) \ar=\ar x(0) + \int_0^t\sigma(x(s)) dB(s) +
\int_0^t\int_0^\infty \phi_0(x(s-)) z \tilde{N}_0(ds,dz) \nnm \\
 \ar \ar
+ \int_0^t b(x(s))ds + \int_0^t\int_0^\infty \phi_1(x(s-))z
N_1(ds,dz) \nnm \\
 \ar \ar
+ \int_0^t\int_{[-1,0)} x(s-)z \tilde{N}_0(ds,dz) + \int_0^t
\int_{[-1,0)} x(s-)z N_1(ds,dz),
 \eeqnn
which is just another form of (\ref{6.3}). \qed

Suppose that $\sigma\ge 0$ and $b$ are real constants. By
Theorem~\ref{t6.3} there is a unique non-negative strong solution
$\{S(t)\}$ to the stochastic differential equation
 \beqnn
dS(t) = \sigma S(t)dB(t) + bS(t)dt + S(t-)dz_0(t) - S(t-)dy_0(t).
 \eeqnn
The process $\{S(t)\}$ is a generalization of the geometric Brownian
motion and has been used widely in mathematical finance; see, e.g.,
Lamberton and Lapeyre (1996, p.144).

Instead of the Lipschitz condition, we may also assume the following
non-Lipschitz condition on the ingredients of (\ref{6.3}):
 \begin{itemize}

\item[{\rm(6.e)}]
The function $x\mapsto \phi_0(x)$ is non-decreasing and there is a
constant $K\ge 0$ so that
 \beqnn
|\sigma(x)-\sigma(y)|^2 + |\phi_0(x)-\phi_0(y)|^2 + |b(x)-b(y)| +
|\phi_1(x)-\phi_1(y)| \leq K|x-y|
 \eeqnn
for all $x,y\ge 0$.

 \end{itemize}

\btheorem\label{t6.4} Under condition (6.e), there is a unique
non-negative strong solution to
 \beqnn
dx(t)
 \ar=\ar
\sigma(x(t))dB(t) + \phi_0(x(t-))dz_0(t) + b(x(t))dt  \nnm \\
 \ar \ar
+\, \phi_1(x(t-))dz_1(t) - x(t-)dy_1(t).
 \eeqnn
\etheorem

We omit the proof of the above theorem since it is a simple
modification of that of Theorem~\ref{t6.3}. A simple generalization
of (\ref{1.8}) is
 \beqnn
dx(t) \ar=\ar \sqrt{2ax(t)} dB(t) + \sqrt[\alpha]{cx(t-)} dz_0(t) +
(\beta x(t)+b)dt \\
 \ar \ar
+\, dz_1(t) - x(t-)dy_1(t).
 \eeqnn
By Theorem~\ref{t6.4} there is a unique non-negative strong solution
to the above equation. It is natural to call the solution a
continuous state process with immigration and emigration
(CBIE-processes). The last term on the right hand side of the
equation brings negative jumps representing the occurrences of
emigration.

\bigskip

\textbf{Acknowledgements.} We would like to thank Professor
L.~Mytnik for pointing out the mistakes in the proofs in an earlier
version of the paper. We are grateful to Professors J.~Jacod and
R.~Schilling for helpful comments on the literature.

%%%%%%%%%%%%%%%%%%%%%%%%%%%%%%%%%%%%%%%

\end{document}